\documentclass{article}
\usepackage{amsfonts}
\usepackage{amsmath, amssymb}
\usepackage{color}
\newtheorem{Theorem}{Theorem}[section]
\newtheorem{Lemma}[Theorem]{Lemma}
\newtheorem{Proposition}[Theorem]{Proposition}
\newtheorem{Definition}[Theorem]{Definition}
\newtheorem{Remark}[Theorem]{Remark}
\newtheorem{Corollary}[Theorem]{Corollary}
\newtheorem{Example}[Theorem]{Example}

\def\RR{{\mathbb R}}
\def\PP{{\mathbf P}}
\def\DD{{\mathbb D}}
\def\dd{{\mathrm d}}
\def\DDD{{\mathrm D}}

\def\o{\omega}
\def\O{\Omega}
\def\va{\varphi}
\def\la{\lambda}
\def\vth{\vartheta}
\def\eps{\varepsilon}

\def\F{{\cal F}}

\def\sD{{\cal D}}

\newcommand{\N}{{\mathbb N}}

\begin{document}

\title{Invariance and Monotonicity for
Stochastic Delay Differential
Equations}

\author{Igor Chueshov\footnote{
Department of Mechanics and Mathematics,
 Kharkov National University, 61077 Kharkov, Ukraine;
e-mail: chueshov@univer.kharkov.ua} $\ $ and Michael Scheutzow\footnote{
Institut f\"ur Mathematik,
Technische Universit\"at Berlin,
Str. des 17 Juni 136, 10623 Berlin, Germany;
e-mail: ms@math.tu-berlin.de}
}

\maketitle

\begin{abstract}
We study invariance and monotonicity properties  of Kunita-type
sto\-chastic differential equations in $\RR^d$ with delay.
Our first result provides sufficient conditions for the invariance
of closed subsets of $\RR^d$.
Then we present a comparison principle and show that under appropriate conditions
the stochastic delay system considered generates a monotone
(order-preserving) random dynamical system. Several applications
are considered.
\end{abstract}

\bigskip


\section{Introduction}
In this paper, we study invariance and monotonicity properties of
 a class of stochastic functional differential equations
(sfde's) driven by a {\em Kunita-type} martingale field.
Our main results are Theorem~\ref{th-inv} on deterministic invariant domains and the comparison principle stated
in Theorem~\ref{th-cp1}. To prove them we represent the
sfde as a random fde (see \cite{mohascheu03} and the references therein).
From the point of view of deterministic delay systems this random
fde has a nonstandard structure and therefore we cannot apply the results on monotonicity available in the deterministic theory. This is
 why we are forced to develop a new method
starting from the basic monotonicity ideas.
We restrict our attention to a class of sfde's which generate
a stochastic semi-flow on the state space of continuous functions
(for an example of an sfde which does not generate such a
semi-flow, see \cite{moha86}). For other classes (on $L_p$-type spaces,  for instance)
we can use a variety of approximation procedures to achieve similar results.
Our choice of continuous functions as a phase space is mainly motivated
by the fact that some  important results
in the theory of monotone systems require a phase space with a
solid minihedral cone (see, e.g., \cite{Kr1,Kr4}).
\par
We note that
invariance properties for deterministic functional differential equations  have been discussed
by many authors (see, e.g., \cite{seifert, martin1, martin2}
and the references therein).
We also refer to \cite{Sm} and to the literature quoted there
 for monotonicity properties of  deterministic fde's.
 Stochastic and random ode's were considered in \cite{Chu02}.
Similar questions for  nonlinear stochastic  partial differential equations
(spde's) were
studied in \cite{Chu00,ChVui2004} (see also
\cite{ChVui1998,ChVui2000} and the references therein
for other applications  of monotonicity methods in spde's).
\par
The paper is organized as follows.
\par
In Section~\ref{sect2} we introduce basic definitions and hypotheses and describe the structure of our
stochastic fde model in \eqref{eq-d}
and its random representation (see \eqref{int-d1}).
The central result in this section is Proposition~\ref{pr:equiv-d}
which shows the equivalence of the stochastic fde \eqref{eq-d}
and the random fde \eqref{int-d1}.

\par
In Section~\ref{sect3} we establish our main result concerning the invariance
of deterministic domains (see Theorem~\ref{th-inv}).
The proof involves the random representation established
in Proposition~\ref{pr:equiv-d} and also the deterministic approach developed in
\cite{seifert}. As an application of Theorem~\ref{th-inv}
we consider an invariant regular simplex for stochastic delayed
Lotka-Volterra type model.
\par
In Section~\ref{sect4} we consider  quasi-monotone vector (drift) fields and using
 the same idea as in Section~\ref{sect3}  establish in Theorem~\ref{th-cp1}
a comparison principle for the corresponding sfde's.
\par
In Section~\ref{sect5}
we apply the results of Sections~\ref{sect3} and \ref{sect4} to construct
{\em random dynamical systems}
(RDS's) defined on invariant regions and
 generated  by sfde's from the class considered (see Theorem~\ref{th-gen}).
 These RDS's become order-preserving for quasi-monotone drift fields (see Theorem~\ref{th-mrs}).
 In this section following \cite{Arn98} (for the monotone
case, see also \cite{Chu02}) we  recall  well-known notions of the theory of random dynamical systems including that of a {\em pull-back attractor}.
Theorem~\ref{th-mrs} on the generation of a monotone RDS allows us apply results from the theory of monotone RDS's
(see, e.g., \cite{ArChu98a,Chu02,ChuSch03a} and the literature cited in these publications)
to describe the qualitative dynamics of the sfde's  considered.
We discuss this issue briefly and provide several examples.

\section{Preliminaries}\label{sect2}

Let $r>0$, $d$ a positive integer and let $C:=C([-r,0],\RR^d)$ be the Banach space of
continuous $\RR^d$-valued functions equipped with the supremum norm $\Vert \cdot \Vert_C$.
For a continuous $\RR^d$-valued function $x$ defined on some subset of $\RR$ containing
the interval $[s-r,s]$, we define $x_s \in C$ by
$$
x_s (u):=x(s+u),\;u \in [-r,0].
$$
Let $(\Omega,{\cal F},\left(\F_t\right)_{t\ge 0},\PP)$ be a filtered
probability space satisfying the usual conditions.
On this probability space we define real-valued random fields $M^i$ and $G^i$, $i=1,2,\dots,d$
satisfying the following hypotheses.\\

\noindent {\bf Hypothesis (M).}  For each $i=1,2,\dots,d$,
$M^i: [0,\infty) \times \RR^d \times \O \to \RR$ satisfies
\begin{itemize}
\item[(i)] $M^i$ is continuous in the first two variables for each $\o \in \O$.
\item[(ii)]  For each $x \in \RR^d$, $M^i(.,x)$ is a local martingale
and  $M^i(0,x,\o)=0$ for all $\o\in\O$.
\item[(iii)] There exist  $\delta \in (0,1)$ and
predictable processes
$a^{ij}: [0,\infty) \times \RR^{2d} \times \O \to \RR$ such that for each $i,j \in \{1,...,d\}:$
\begin{multline*}
 R^{ij}(t,\o):=
 \sup_{x,y \in \RR^d} \frac{\vert a^{ij}(t,x,y)\vert}{(1+|x|)(1+|y|)}
+ \sup_{x,y \in \RR^d}
\| \DDD_x \DDD_y a^{ij}(t,x,y) \|  \\
  +\;
\sup_{x\neq x',y \neq y'} \frac{\|\bar a^{ij}(t,x,y)-\bar a^{ij}(t,x,y')-
\bar a^{ij}(t,x',y)+\bar a^{ij}(t,x',y')\|}{|x-x'|^{\delta}|y-y'|^{\delta}}
\end{multline*}
is finite, where $\bar a^{ij}(t,x',y'):=\DDD_x \DDD_y a^{ij}(t,x',y')$ and
\[
\langle M^i(\cdot,x),M^j(\cdot,y) \rangle_t=\int_0^t a^{ij}(s,x,y,\o)\,\dd s
~~{\rm
a.s.,}~~i,j=1,\ldots, d,
\] where $\langle M^i,M^j \rangle_t$ denotes the corresponding
joint quadratic variation (see \cite{kunita} for details).
Moreover, we assume that
 the map $t \mapsto R^{ij}(t,\o)$
is locally integrable w.r.t. Lebesgue measure
for every $\o \in \O$ and $i,j \in \{1,...,d\}$. In the definition of $R^{ij}$,
$D_xD_y a$ denotes the matrix formed by the corresponding partial derivatives and
$\|.\|$ is an arbitrary norm on the space of
matrices.
\end{itemize}

\noindent {\bf Hypothesis (G).} $G=(G^1,...,G^d): [0,\infty)\times C \times \O \to \RR^d$ satisfies
\begin{itemize}
\item[(i)] $G^i$ is jointly continuous in the first two variables for each $\o \in \O$.
\item[(ii)] For each $\o \in \O$, bounded set $B$ in $C$  and $T>0$ there exists some $L=L(T,B,\o)<\infty$ such that
$\vert G^i (t,\eta,\o)- G^i (t,\zeta,\o) \vert \le L\Vert \eta-\zeta \Vert_C$ for all
$0 \le t \le T$ and $\eta,\zeta \in B$.
\item[(iii)] For each $\eta \in C$ and $t \in [0,\infty)$, $G (t,\eta)$ is $\F_t$-measurable.
\end{itemize}
Below, it will be important to decompose $G$ as
\[
G(t,\eta,\o)=H(t,\eta,\o)+b(t,\eta(0),\o),
 \]\
where both $H$ and $b$ satisfy (i), (ii) and (iii)
of the previous hypothesis (with $C$ replaced by $\RR^d$ with the Euclidean norm for $b$). In addition, we assume that
$b(t,.)$ is continuously differentiable for each $t$ and $\o$ and there exist $\delta >0$ and a number
$c(T,\o)< \infty$ such that
\begin{equation}\label{b}
\sup_{0 \le t \le T}\Big\{\sup_{x \in \RR^d}\|\DDD b(t,x,\o)\| + \sup_{x,y \in \RR^d,x \neq y} \frac {\|\DDD b(t,x,\o)-\DDD b(t,y,\o)\|}{|x-y|^{\delta}}\Big\} \le c(T,\o).
\end{equation}
In this case we say that Hypothesis {\bf (G)} holds with decomposition $G=H+b$.

For the rest of this section, we assume both hypotheses {\bf (M)} and {\bf (G)} and fix a particular decomposition $G=H+b$ as above.

We consider the following Kunita-type delay stochastic differential equation
\begin{align}\label{eq-d}
\left\{
\begin{array}{rll}
\dd x^i(t)&=&G^i(t,x_t)\,\dd t + M^i(\dd t,x(t)),\quad i=1,2,\dots,d, \quad t\ge s,\\
x_s&=&\eta,
\end{array}
\right.
\end{align}
where $s \ge 0$ and $\eta$ is a $C$-valued $\F_s$-measurable random variable.

For the definition of {\em Kunita-type stochastic integrals}
$$
\int_s^t M^i(\dd u, x(u)),
$$
for adapted and continuous (or more general) processes $x$, the reader is referred to Kunita's monograph \cite{kunita}.
Readers who are unwilling to learn Kunita integrals (even though they are very natural and easy to deal with objects) can
think of the special case
\begin{equation}\label{finite}
M^i(t,x):=\sum_{k=1}^m \int_0^t \sigma^{ik}(s,x)    \dd W^k(s),
\end{equation}
where $W^k,\,k=1,...m$, are independent Brownian motions and the $\sigma^{ik}$ are (deterministic) functions (satisfying appropriate
regularity properties). In this case \eqref{eq-d} reads
\begin{align*}
\left\{
\begin{array}{rll}
\dd x^i(t)&=&G^i(t,x_t)\,\dd t + \displaystyle\sum\limits_{k=1}^m \sigma^{ik}(t,x(t))\dd W^k(t),~~~ i=1,\ldots,d, ~~ t\ge s,\\ [2 mm]
x_s&=&\eta,
\end{array}
\right.
\end{align*}
and $a^{ij}(t,x,y)=\sum_{k=1}^m \sigma^{ik}(t,x)\sigma^{jk}(t,y)$ is deterministic.

We aim at a representation of the solution from which one can read off continuity properties with respect to the initial condition. Note that
even though equation \eqref{eq-d} is easily seen to have a unique solution for each fixed $s$ and $\eta$, continuity with respect to
$\eta$ does not follow since solutions are defined only up to a set of measure zero which may depend on $\eta$. To obtain continuity, one has to
select a particular modification of the solution. We will use a variant of the {\em variation-of-constants} technique which turns \eqref{eq-d}
into an equation which does not contain any stochastic integral and can therefore be solved for each fixed $\o \in \O$. We will see that the
modification of the solution which is given by the pathwise equation does automatically exhibit continuous dependence upon the initial condition.
The  variation-of-constants technique, which  is well-known for ode's, has already been applied to sfde's in \cite{moha90} and \cite{mohascheu03}.

For further use we need some properties of the following (non-delay) stochastic  equation
\begin{align}
\label{psi-d}
\left\{
\begin{array}{rll}
\dd \psi^i(t)&=&b^i(t,\psi(t))\,\dd t + M^i(\dd t,\psi(t)),\quad i=1,2,\dots,d, \quad t\ge s,\\
\psi^i(s)&=&x,
\end{array}
\right.
\end{align}
where $\psi=(\psi^1,\ldots,\psi^d)$.
The following lemma states that equation \eqref{psi-d} generates a stochastic flow of
diffeomorphisms in $\RR^d$.
This is a special case of Theorem 4.6.5 in \cite{kunita}.

\begin{Lemma}\label{le:diff}
We assume that $b\equiv(b^1,\ldots,b^d)\, : [0,\infty) \times \RR^d \times \Omega \to \RR^d$ is a vector field satisfying \eqref{b}.
Then there exists a process
$\Psi: [0,\infty)^2 \times \RR^d \times \O \to \RR$ which satisfies the following:
\begin{itemize}
\item[(i)] For each $s \ge 0$ and $x \in \RR^d$,
$\psi(t)\equiv \Psi_{s,t}(x,\o)$, $t \ge s$ solves equation \eqref{psi-d}.
\item[(ii)] For each $s \ge 0$, $x\in\RR^d$ and $\o \in \O$,
$\Psi_{s,s}(x,\o)=x$.
\item[(iii)] The maps $(s,t,x) \mapsto
\Psi_{s,t}(x,\o)$ and $(s,t,x) \mapsto \DDD_x\Psi_{s,t}(x,\o)$
are continuous for each $\o \in \O$. Furthermore
$\Psi_{s,t}(.,\o)$ is a $C^1$-diffeomorphism for each $s,t \in \RR$
and  $\o \in \O$.
\item[(iv)] For each $s,t,u \ge 0$, and
$\o \in \O$ we have the semi-flow property
\[
\Psi_{s,u}(\cdot,\o)=\Psi_{t,u}(\cdot,\o) \circ \Psi_{s,t}(\cdot,\o).
\]
\end{itemize}
\end{Lemma}
Note that by (ii) and (iv), we have $\Psi_{s,t}(.\o)=\big( \Psi_{t,s}(.,\o) \big)^{-1}$.
\smallskip\par
Lemma~\ref{le:diff}  allows us to construct the following representation
for  solutions to \eqref{eq-d}. In the special case in which the martingale field $M$ is given by a finite number of
Brownian motions as in \eqref{finite} and $b \equiv 0$, this representation was established in Lemma 2.3 in \cite{mohascheu03}.

Let $\Psi(u,x,\o):=\Psi_{0,u}(x,\o)$. We define the functions
$\xi:[0,\infty) \times \RR^d \times \O \to \RR^d$ and
$F: [0,\infty) \times \RR^d \times C \times \O \to \RR$ by
\begin{eqnarray*}
\xi(u,x,\o)&:=&\Psi (u,\cdot,\o)^{-1}(x)=\Psi_{u,0}(x,\o),\\
F(u,x,\eta,\o)&:=&\left\{ \DDD_x\Psi(u,x,\o) \right\}^{-1} H(u,\eta,\o)
\end{eqnarray*}
and consider the (random) equation
\begin{align}
\label{int-d1}
x(t,\o)=
\Psi\Big( t,\Big[ \xi(s,\eta (0),\o)+\int_s^t F(u,\xi (u,x(u,\o),\o),x_u(\o),\o)\dd u
\Big],\o \Big)
\end{align}
for $t \ge s$ with the initial data
\begin{align}\label{int-d2}
x(t,\o)=\eta(t-s) \quad\mbox{for}\quad t \in [s-r,s],
\end{align}
where $\eta$ is a $C$-valued $\F_s$-measurable random variable.
We suppress the dependence of $x=(x^1,\dots,x^d)$ on $s$ and $\eta$ for notational simplicity.
The following proposition shows that equations \eqref{eq-d}
and \eqref{int-d1} (together with \eqref{int-d2}) are equivalent.
\begin{Proposition}\label{pr:equiv-d}
Fix $s \ge 0$, a $C$-valued $\F_s$-measurable random variable $\eta$
and a stopping time $T \ge s$.
An adapted $\RR^d$-valued process $x(t)$ with
continuous paths solves equation \eqref{eq-d}
on the interval $[s,T(\o)]\cap [s,\infty)$ with
initial condition $x_s=\eta$ if and only if  $x$
satisfies  \eqref{int-d1} and  \eqref{int-d2} on the same
interval for  almost all $\o \in \O$.
\end{Proposition}

\noindent {\bf Proof.}
The proof is essentially the same as that of Lemma 2.3 in \cite{mohascheu03}. Our assumptions are slightly different from the ones in
\cite{mohascheu03} but this does not affect the arguments in the proof. Therefore, we skip some details.

First assume that $x$ solves  \eqref{int-d1} and  \eqref{int-d2} on
$[s,T(\o)]\cap [s,\infty)$ for almost all $\o \in \O$. Then $x_s=\eta$ almost surely. Equation \eqref{int-d1},
together with a slight modification of the generalized It\^o's formula as stated in \cite{kunita}, Theorems 3.3.1 and 3.3.3(i)
imply, that $x$ is a continuous semimartingale and satisfies
\begin{align*}
\dd x(t)=& \DDD_x\Psi\left(t,\big[ \xi(s,\eta (0),\o)+\int_s^t F(u,\xi (u,x(u,\o),\o),x_u(\o),\o)\dd u
\Big],\o \right)\\
&\times F \left(t,\xi(t,x(t),\o),x_t,\o\right)\dd t\\
&+\Psi\big( \dd t,\big[ \xi(s,\eta (0),\o)+\int_s^t F(u,\xi (u,x(u,\o),\o),x_u(\o),\o)\dd u
\Big],\o \big)\\
=&H(t,x_t)\,\dd t + b(t,x(t))\dd t + M(\dd t,x(t)).
\end{align*}
Therefore $x$ solves \eqref{eq-d}.
\par
Conversely, suppose that $x$ solves \eqref{eq-d} and define
$$
\zeta(t,\o):=\xi(s,\eta(0),\o)+\int_s^t F(u,\xi(u,x(u),\o),x_u,\o)\dd u.
$$
Let
\begin{align*}
\tilde x(t,\o)\equiv (\tilde x^1(t,\o),\ldots, \tilde x^d(t,\o)) :=\left\{
\begin{array}{ll}
\Psi \left( t,\zeta(t,\o),\o\right), \quad &t \ge s\\
\eta (t-s), \quad &t \in [s-r,s].
\end{array}
\right.
\end{align*}
One can see that $\tilde x^i(t,\o)$
is a semimartingale with differential
$$
\dd\tilde x^i(t)=H^i(t,x_t)\dd t + b^i(t,\tilde x(t))\dd t + M^i(\dd t,\tilde x(t)), \quad i=1,\ldots, d.
$$
This (non-retarded) sde has a unique solution $\tilde x$
with initial condition
$ \tilde x(s)=\eta(0)$, so   $x$ and $\tilde x$ agree on $[s-r,T] \cap [s,\infty)$ almost surely.
This proves the proposition. \hfill $\Box$
\smallskip\par
The following  proposition provides a well-posedness result
concerning problem \eqref{int-d1} and \eqref{int-d2}.

\begin{Proposition}\label{pr:wp}
Let hypotheses $\mathrm{({\bf M})}$ and $\mathrm{({\bf G})}$ be satisfied.
Then there exists a set $\Omega_0$ of full measure such that
for all $\omega \in \Omega_0$, $s \ge 0$, and $\eta \in C$,
the problem \eqref{int-d1} and \eqref{int-d2} has a unique local
solution $x(s,\eta,t,\o)$ up to an explosion time $\tau(s,\eta,\o)$. For each $s \ge 0$, the solution depends continuously upon $(t,\eta)$ (up to explosion).
Further, the following semi-flow property holds:
for all $0 \le s \le t \le u$, all $\eta \in C$ and all $\o \in \O_0$, we have
$$
x(s,\eta,u,\o)=x(t,x_t(s,\eta,t,\o),u,\o)~~~\mbox{up to explosion.}
$$

\end{Proposition}

\noindent {\bf Proof.} This is (essentially) Theorem 2.1 in \cite{mohascheu03}. The only differences
are the fact that in \cite{mohascheu03} the authors use the Hilbert space M$_2$ instead of $C$ as the state space and that we separate
$b$ from $G$ and combine it with the martingale part. The proof of local existence, uniqueness and continuity of the problem
\eqref{int-d1} and \eqref{int-d2} is based on a rather standard fixed point argument.
\hfill $\Box$\\

It is natural to ask for sufficient conditions for the explosion time $\tau(s,\eta,\o)$ to be infinite on a set of full measure which does not depend on
$s$ and $\eta$.
We will say that condition $\mathrm{({\bf GE})}$ (for {\em global existence}) holds if $\mathrm{({\bf G})}$ and $\mathrm{({\bf M})}$ hold with decomposition
$G=H+b$ and there exists a
set $\Omega_0$ of full measure such that $\tau(s,\eta,\o)=\infty$ for all $s \ge 0$, all $\eta \in C$ and all $\o \in \Omega_0$.
Various sufficient conditions for  $\mathrm{({\bf GE})}$ are formulated in
Theorem 3.1 in \cite{mohascheu03}. They are based on spatial estimates  on the growth of the flow $\Psi$ and its spatial derivative
which were established in \cite{mohascheu98} and \cite{imscheu99}. We quote them here:

\begin{Proposition}\label{criteriaGE}
Let  $\mathrm{({\bf G})}$ and $\mathrm{({\bf M})}$ hold with decomposition $G=H+b$. Each of the following conditions is sufficient for $\mathrm{({\bf GE})}$:
\begin{itemize}
\item[(i)] For each $T>0$ and $\o \in \O$ there exist $c=c(T,\o)$ and
$\gamma=\gamma(T,\o) \in [0,1)$ such that
\begin{equation}\label{unno}
|H(t,\eta,\omega)| \le c(1+\|\eta\|_C^{\gamma})
\end{equation}
for all $0 \le t \le T$, $\eta \in C$ and $\omega \in \Omega$.
\item[(ii)] For each  $T>0$ there exists $\beta \in (0,r)$ such that $H(u,\eta,\o)=H(u,\tilde \eta,\o)$ holds for all $\o \in \Omega$ whenever
$0 \le u \le T$ and $\eta |_{[-r,-\beta]}=\tilde  \eta |_{[-r,-\beta]}$.
\item[(iii)]  For all $\o \in \Omega$ and $T \in (0,\infty)$ we have that
 \[
\sup_{0 \le u \le T,x \in \RR^d}\|(\DDD_x\psi(u,x,\o))^{-1}\|<\infty
 \]
 and there exists $c=c(T,\o)$
such that \eqref{unno} holds with $\gamma=1$.
\end{itemize}
\end{Proposition}

It is a bit annoying that (i) excludes the case of  $H$ satisfying a global Lipschitz condition. It seems to be open whether  $\mathrm{({\bf GE})}$ holds in that case.

\section{Deterministic invariant regions}\label{sect3}
In this section we assume that Hypotheses  ({\bf M}) and ({\bf G}) with decomposition $G=H+b$
and condition \eqref{unno} are in force and consider a general problem of the form \eqref{int-d1},\eqref{int-d2}.
We provide sufficient conditions that, given a non-empty closed (deterministic) subset $\DD$  in $\RR^d$, a solution
with values in this set for $t\in [t_0-r,t_0]$ will have values in $\DD$ for all $t>t_0$.
The key idea is to decompose the solution semi-flow in such a way that $\Psi$ alone leaves $\DD$ invariant and that the remaining drift does not
change this property.
\par
Below we use the notation
\begin{equation}
\label{d-0}
C_\DD=\left\{ \eta\in C\, :\,  \eta(s)\in\DD~~~\mbox{for every}~
s\in  [-r,0]\right\}.
\end{equation}
We need some additional hypotheses (which are inspired by similar hypotheses for deterministic fde's in \cite{seifert}).\\

\noindent {\bf Hypothesis (G$_\eps$).}
There exists a family $\{ G_\eps\}$ of random fields satisfying {\bf (G)}
with decomposition $G_\eps=H_\eps+b$ for every $\eps\in (0,\eps_0]$ such that
\begin{itemize}
\item[(i)] $H_\eps$ satisfies condition \eqref{unno};
\item[(ii)] $\lim_{\eps\to 0} H_\eps(t,\eta,\o)= H(t,\eta,\o)$
for every $(t,\eta)\in [0,\infty) \times C_\DD$ and $\o\in\O$;
\item[(iii)]
given $(t,\eta)\in [0,\infty) \times C_\DD$
 and $\eps\in (0,\eps_0]$ there exists
an $\alpha=\alpha(\eps,t,\eta,\o)>0$ such that if $0<h\le\alpha$
and $u\in\RR^d$ is such that $|u|\le\alpha$, then
\begin{equation*}
\eta(0)+hH_\eps(t,\eta,\o)+hu\in\DD;
\end{equation*}
\item[(iv)] if $y^\eps(t,\eta)$ solves the problem
\eqref{int-d1},\eqref{int-d2} with $G_\eps$ instead of $G$,
then  for every $\o\in\O$, $t \ge 0$, and $\eta \in C$ we have
$\lim_{\eps\to 0} y^\eps(t,\eta,\o)= x(t,\eta,\o)$,
where  $x(t,\eta,\o)$ solves \eqref{int-d1},\eqref{int-d2}.
\end{itemize}

 Note that condition (ii) implies that $\DD$ is the closure of its interior.
Further note that condition (iv) above is not implied by the other conditions -- not even in the case of a deterministic ode, see \cite{seifert}.

\noindent {\bf Hypothesis (M$_{\DD}$).}  The problem \eqref{psi-d}
generates a stochastic flow $\Psi_{t,s}(\cdot,\o)$ of
diffeomorphisms of $\RR^d$ such that
\begin{equation}\label{inv}
\Psi_{s,t}(\DD,\o)=\DD,~~ t>s,\; \o\in\O.
\end{equation}

\begin{Remark}\label{kuni}
{\rm For flows which are driven by a finite number of Brownian motions, explicit criteria for the validity of
this hypothesis are well-known (we will state some of them below). We have not found corresponding criteria for Kunita-type equations
in the mathematical literature. In fact, such criteria follow easily in case $\DD$ is compact: for Kunita-type sde's, the one-point motion (i.e.~the solution for
a single starting point $x \in \RR^d$) can be described by an equivalent sde which is driven by a finite number of Brownian motions (which depend on
the point $x$). Assuming that for each $x \in \DD$ the solution stays in $\DD$ forever with probability one (for which one can check the
known criteria), then the same holds true for a countable dense set of initial conditions in $\DD$. The fact that $\DD$ is compact and the flow is
continuous shows that there exists a set $\Omega_0$ of full measure such that \eqref{inv} holds for all $\omega \in \Omega$. Our claim follows
since we are free to modify
$\Psi$ on a set of measure zero.}
\end{Remark}


\begin{Remark} \label{re3.2}
{\rm
If $\DD$ is a closed convex subset of $\RR^d$ with nonempty interior
then  Hypothesis {\bf (G$_\eps$)} follows from the Nagumo type relation
\begin{equation}\label{d-2}
\lim_{h\to 0+} h^{-1}{\rm dist}(\eta(0)+hH(t,\eta,\o),\DD)=0
\end{equation}
for any $t\in \RR$, $\o\in\O$ and $\eta\in C$ such that $\eta(s)\in\DD$
for $s\in [-r,0]$. In this case we can take
\[
H_\eps(t,\eta,\o)=H(t,\eta,\o)-\eps (\eta(0)-{\bf e}),
\]
where ${\bf e}$ is an element from int$\,\DD$.
If $\DD=\RR^d_+$ relation \eqref{d-2}
is equivalent to the requirement
\begin{equation}\label{d-3}
\{\eta\ge 0,\; \eta^i(0)=0\}\;
\Rightarrow\;
H^i(t,\eta,\o)\ge 0
\end{equation}
for every $t$, $\o\in\O$ and $i=1,\ldots,d$.
For the proofs we refer to \cite{seifert}.
}
\end{Remark}

In the following remark  we discuss conditions and examples when Hypothesis
{\bf (M$_{\DD}$)}  is valid.

\begin{Remark}\label{re3.1}
{\rm
Assume that $\DD$ is a closed set in $\RR^d$ such that $\DD$
has an outer normal at every point of its boundary.
 We recall that
a unit vector $\nu$ is said to be an {\em outer normal} to $\DD$ at the point
$x_0\in \partial \DD$, if there exists a ball
$B(x_1)$ with  center at $x_1$ such that
$B(x_1)\cap \DD=\{ x_0\}$ and
$\nu=\lambda\cdot (x_1-x_0)$ for some positive
$\lambda$.
\par
Let $W_1,...,W_l$ be independent standard Wiener processes.
We consider \eqref{psi-d}   with
 \begin{equation}\label{d-1a}
M^i(\dd t,\psi(t))=\sum_{j=1}^lm^i_j(\psi(t))  \dd W_j(t),
\end{equation}
where the coefficients  have bounded derivatives up to second order.
The problem in \eqref{psi-d} can be written as
a Stratonovich sde:
\begin{align*}
\left\{
\begin{array}{rll}
\dd \psi^i(t)&=&\tilde{b}^i(t,\psi(t))\,\dd t + \displaystyle\sum\limits_{j=1}^lm^i_j(\psi(t)) \circ \dd W_j(t),~~ i=1,\dots,d,
\\ [2mm]
\psi^i(s)&=&x,
\end{array}
\right.
\end{align*}
where ``$\circ$'' denotes Stratonovich integration
and
\[
\tilde{b}^i(t,x)\equiv b^i(t,x) -\frac12\sum_{j=1}^l\sum_{k=1}^d
 m^k_j(x) \frac{\partial m^i_j(x)}{\partial x_k}.
\]
It follows from Wong-Zakai type arguments that $\DD$ is forward invariant under $\psi$
if for any $x\in \partial \DD$ we have
\begin{equation}\label{1.9.4}
\sum_{i=1}^d \tilde{b}^i(t,x) \nu^i_{x}
\le 0\quad\mbox{and}\quad
\sum_{i=1}^d m^i_j(x) \nu^i_{x} = 0,\quad j=1,\ldots,l,
\end{equation}
for every outer normal $\nu_{x}=(\nu^1_{x},\ldots,\nu^d_{x})$ to $\DD$
at $x$.   We refer to \cite[Chap.2, Corollaries 2.5.1 and 2.5.2 ]{Chu02} for details.
Further, {\bf (M$_{\DD}$)} holds if $\DD$ is both forward and backward invariant under $\psi$. For this
to hold it is sufficient to assume that the first inequality in \eqref{1.9.4} is an equality for each
$x \in \partial \DD$.
We note that in the case $\DD=\RR_+^d$ and  $b^i(t,x)\equiv0$ the first condition
in \eqref{1.9.4} follows from the second one
which can be written in the form
\[
m^i_j(x)=0~~{\rm for~all}~~ x=(x_1,\ldots,x_{i-1},0,x_{i+1}\ldots,x_d),~~
i=1,\ldots, d,~j=1,\ldots,l.
\]
As an example we point out the case when
$\DD=\{ (x_1;x_2): x^2_1+x^2_2\le 1\}\subset\RR^2$ and
problem \eqref{psi-d} has the form
\[
\dd x_i= m^i(x_1,x_2) \dd W(t),\quad i=1,2.
\]
In this case $b^i(t,x)\equiv0$  and relations \eqref{1.9.4} holds if $m^1(x) =m^2(x)=0$
for all $|x|=1$.
For instance, we can take
 $m^1(x_1,x_2)=-c(|x|)x_2$ and
$m^2(x_1,x_2)=c(|x|)x_1$, where $c(r)=0$ for  $r=1$.
}
\end{Remark}

Our main result in this section is the following theorem.
\begin{Theorem}\label{th-inv}
Assume that Hypotheses  {\bf (G)}, {\bf (M)},  {\bf (M$_{\DD}$)} and {\bf (G$_\eps$)} hold.
Then $\DD$ is a forward invariant set for  problem \eqref{eq-d} in the sense  that
for any $\o \in \O$, $s \ge 0$ and $\eta\in C$ such that  $\eta(u)\in\DD$
for $u\in [-r,0]$, the (unique) solution $x$ of \eqref{int-d1} and \eqref{int-d2}
satisfies $ x(t,\eta,\o)\in\DD$ for all $t\ge s$.
\end{Theorem}
\par\noindent
{\bf Proof.}
Assume that $\DD$ is not forward invariant.
Then there exist $\o \in \O$, $s \ge 0$, $\eta\in C_\DD$ and $t_*>s$ such that
$x(t_*)\not\in\DD$.
Let $y^\eps(t)$ be a solution to the auxiliary problem in  ({\bf G}$_\eps$)(iv).
It follows from assumption (iv) in ({\bf G}$_\eps$) that there exist
$\eps>0$ and $t_0\in [s,t_*)$ such that
\[
y^\eps(t)\in\DD,\; t\in [s-r,t_0]~~\mbox{and}~~
y^\eps(t_0+h_j)\not\in\DD,
\]
where $\{h_j\}$ is a sequence of positive numbers such that
$\lim_{j\to\infty}h_j=0$. The solution $y^\eps(t)$ can be represented
in the form
\[
y^\eps(t,\o)=\psi(t,\zeta^\eps(t,\o),\o),~~ t\ge t_0,
\]
where
\[
\zeta^\eps(t,\o)=\xi(t_0,y^\eps(t_0),\o)+
\int^t_{t_0} F_\eps(u,\xi(u,y^\eps(u),\o),y^\eps_u,\o)\dd u.
\]
Here
$\xi(t,x,\o)=\psi(t,\cdot,\o)^{-1} x$,
$\psi(t,\cdot,\o)=\Psi_{0,t}(\cdot,\o)$, where
$\Psi_{0,t}(\cdot,\o)$ is the diffeomorphism given by \eqref{psi-d}
and
\[
F_\eps(u,x,\eta,\o)=\left\{ \DDD_x\psi(u,x,\o) \right\}^{-1} H_\eps(u,\eta,\o).
\]
Since $(u,x,\eta) \mapsto F_\eps(u,x,\eta,\o)$ is continuous,
we have that
\[
\int^{t_0+h}_{t_0} F_\eps(u,\xi(u,y^\eps(u),\o),y^\eps_u,\o)\dd u=
h  F_\eps(t_0,\xi(t_0,y^\eps(t_0),\o),y^\eps_{t_0},\o)+ o(h).
\]
Thus
\[
\zeta^\eps(t_0+h,\o)=\xi(t_0,y^\eps(t_0),\o)+
h  F_\eps(t_0,\xi(t_0,y^\eps(t_0),\o),y^\eps_{t_0},\o)+ o(h).
\]
We have that $\psi(t_0,\psi^{-1}(t_0,x,\o),\o)=x$. Therefore by the
chain rule
\[
\DDD_z\psi(t_0,\psi^{-1}(t_0,x,\o),\o)\DDD_x \psi^{-1}(t_0,x,\o)=\mathrm{Id}.
\]
Thus
\[
\left\{\DDD_z\psi(t_0,\psi^{-1}(t_0,x,\o),\o)\right\}^{-1}=
\DDD_x \psi^{-1}(t_0,x,\o).
\]
Consequently,
\[
 F_\eps(t_0,\xi(t_0,y^\eps(t_0),\o),y^\eps_{t_0},\o)=
 \DDD_x\psi^{-1}(t_0,y^\eps(t_0),\o)
 H_\eps(t_0, y^\eps_{t_0},\o).
\]
It is also clear that
\begin{multline*}
\psi^{-1}(t_0,y^\eps(t_0)+ hH_\eps(t_0, y^\eps_{t_0},\o),\o)  -
\psi^{-1}(t_0,y^\eps(t_0),\o) \\
=
h \DDD_x\psi^{-1}(t_0,y^\eps(t_0),\o)
 H_\eps(t_0, y^\eps_{t_0},\o) +o(h).
\end{multline*}
Thus
\[
\zeta^\eps(t_0+h,\o)=
\psi^{-1}(t_0,y^\eps(t_0)+ hH_\eps(t_0, y^\eps_{t_0},\o),\o)
+o(h).
\]
This implies that
\begin{equation}\label{d-4}
\psi(t_0, \zeta^\eps(t_0+h,\o),\o)=
y^\eps(t_0)+ hH_\eps(t_0, y^\eps_{t_0},\o)
+o(h).
\end{equation}
Hypothesis ({\bf G}$_\eps$)(iii)  implies that the right-hand side of \eqref{d-4} is in $\DD$ for all
sufficiently small $h>0$. By  Hypothesis {\bf (M$_{\DD}$)} we therefore have
$\zeta^\eps(t_0+h,\o) \in \DD$ and hence $y^\eps(t_0+h,\o)=\psi(t_0+h,\zeta^\eps(t_0+h,\o) ,\o) \in \DD$ for all sufficiently small
$h>0$ contradicting our assumption that   $y^\eps(t_0+h_j,\o)\notin \DD$ for all $j$.
This contradiction proves the theorem.
\hfill $\Box$
\smallskip\par
In the following assertion we show that,
similarly to the deterministic situation (see \cite{Sm} and the references therein),
in some cases the Nagumo type condition in \eqref{d-2}
provides us necessary and sufficient conditions for invariance.

\begin{Corollary}\label{co:convex-D}
Let $\DD$ be a closed convex subset of $\RR^d$ with nonempty interior
 and let Hypotheses {\bf (G)}, {\bf (M)} and  $\mathrm{(}${\bf M$_{\DD}$}$\mathrm{)}$ be in force.
Then $\DD$ is a forward invariant set if and only if  \eqref{d-2} holds.
\end{Corollary}
{\bf Proof.}
If  \eqref{d-2} holds, then we can apply Remark \ref{re3.2}
to conclude that $\DD$ is forward invariant.
\par
 Let $\DD$ be forward invariant. It is clear that \eqref{d-4} holds for $\eps=0$ and
$t_0=s$, i.e. we have
\begin{equation}\label{d-5}
\psi(s, \zeta(s+h,\o),\o)=
\eta(0)+ hH(s,\eta,\o)
+o(h)
\end{equation}
for any $\eta\in C_\DD$. Since $\DD$ is invariant,
we have $\psi(s+h, \zeta(s+h,\o),\o)=x(s+h,\eta,\o)\in \DD$ for all $h \ge 0$.
({\bf M$_{\DD}$}) implies that $\psi(s, \zeta(s+h,\o),\o)$ lies in $\DD$ for all $h \ge 0$.
Therefore \eqref{d-5} implies  \eqref{d-2}.
\hfill $\Box$
\smallskip\par
In the case $\DD=\RR_+^d$ Corollary~\ref{co:convex-D}
implies the following assertion.
\begin{Corollary}\label{co3.6}
Let {\bf (G)} and {\bf (M)}
be in force with
$M^i$ of  the form \eqref{d-1a}.
Assume that
\[
b^i(t,x)= 0,~~~
m^i_j(x)=0~~\forall\; x=(x_1,\ldots,x_{i-1},0,x_{i+1}\ldots,x_d),~~
\]
where $i=1,\ldots, d$, $j=1,\ldots,l$.
Then $\RR^d_+$ is forward invariant set if and only if \eqref{d-3} holds.
\end{Corollary}
{\bf Proof.} It follows from Corollary~\ref{co:convex-D}, see also Remarks \ref{re3.1} and \ref{re3.2}.
\hfill $\Box$
\smallskip\par
More complicated example of an invariant set $\DD$ is discussed in the
following remark.

\begin{Remark} \label{re3.3}
{\rm
Assume that $\DD$ is a set of the form
\[
\DD=\left\{ x\in \RR^d\, :\, \langle a_q, x\rangle\le \gamma_q,\;
q=1,\ldots, Q
 \right\},
\]
where $a_q=(a_q^1,\ldots, a_q^d)\in\RR^d$, $\gamma_q\in\RR$,
$q=1,\ldots, Q$. Then \eqref{d-2} holds if and only if
for every $q=1,\ldots, Q$ we have the relation
\[
\sum_{i=1}^d a^i_q H^i(t,\eta,\o)\le 0
\]
whenever $\eta\in C_\DD$ and  $\langle a_q, \eta(0)\rangle= \gamma_q$.
\par
In the case when $b^i\equiv 0$ and $M^i$ has the form \eqref{d-1a},
it follows from Remark~\ref{re3.1}
that Condition {\bf (M${}_\DD$)} holds
if
for every $q=1,\ldots, Q$ we have the relations
\[
\sum_{j=1}^l\sum_{k=1}^d
\sum_{i=1}^d m^k_j(x) \frac{\partial m^i_j(x)}{\partial x_k} a^i_q
= 0,\qquad~~
\sum_{i=1}^d m^i_j(x) a^i_q = 0,\quad j=1,\ldots,l,
\]
whenever $\langle a_q, x\rangle= \gamma_q$.
For instance, this condition is true if
\[
m^i_j(x)=\sigma_j^i(x_1,\ldots,x_d)\cdot
\prod_{q=1}^Q  h_q(\langle a_q, x\rangle- \gamma_q),
\]
where $\sigma_j^i(x_1,\ldots,x_d)$ are arbitrary and
$h_q(s)$ is such that $h_q(0)=0$.
}
\end{Remark}
Now we provide some examples.

\begin{Example}\label{ex3.1}
{\rm
For the system
\[
\dd x^i(t)= x^i(t) f^i(x_t)\,\dd t +x^i(t)\sum_{j=1}^l\sigma_j^i(x^1(t),\ldots,
x^d(t))\,\dd W_j(t),~~ i=1,\ldots, d,
\]
the set $\RR^d_+$ is a forward invariant set.
Here $f^i$ and $\sigma_j^i$ are such that condition {\bf (GE)} holds.
This conclusion follows from Corollary \ref{co3.6}.
}
\end{Example}

\begin{Example}
{\rm
In the previous example, the noise can also be replaced by a more general {\em Kunita-type} noise.
As a particular example, let $N$ be {\em space-time white noise} on $\RR^d \times [0,\infty)$, let
$h:\RR^d \to [0,\infty)$ be $C^{\infty}$ with compact support and define
$$
M^i(\dd t,x):=\phi(x^i) \int_{\RR^d} h(x-z) N(\dd z,\dd t),
$$
where $\phi \in C^{\infty}$ is bounded and all its derivatives are bounded. Assume that $\phi(0)=0$.
Then
$$
\langle M^i(.,x), M^j(.,y)\rangle_t=\phi(x^i) \phi(y^j)  t \int_{\RR^d} h(x-z) h(y-z) \dd z.
$$
Note that $M$ satisfies hypothesis {\bf (M)} due to our assumptions on $\phi$ and $h$.
To facilitate things, we assume that $f^i(\eta)=g^i(\eta(-r))$ where $g^i$ is bounded and Lipschitz.
Then the set $\RR^d_+$ is a forward invariant set for the flow generated by
$$
\dd x^i(t)= \phi(x^i(t)) g^i(x(t-r))\,\dd t + M^i(\dd t,x(t)),\,i=1,...,d.
$$
To see this, first note that condition {\bf (GE)} holds by Proposition
\ref{criteriaGE}(i). Then we argue as in Remark \ref{kuni}:
for each given starting point $x \in \RR^d_+$, the solution starting in $x$ will remain in  $\RR^d_+$ forever almost surely since it
can be written as a solution to an sfde with finitely many driving Brownian motions, so the same property holds for all starting points
in  $\RR^d_+$ with rational coordinates. Since we know that a local flow exists, no trajectory of the flow can leave  $\RR^d_+$
through one of the hyperplanes bordering  $\RR^d_+$. Since also  {\bf (GE)} holds, no trajectory of the
flow can escape to infinity at finite time either, so the set $\RR^d_+$ is invariant.
\par
Note that in this set-up the driving noise $M$ is independent at locations $x$ and $y$ with distance larger than the diameter of
the support of $h$ which is a reasonable assumption in many models and which cannot be achieved with a finite number of
driving Wiener processes.
}
\end{Example}

\begin{Example}[Lotka-Volterra type model]\label{ex3.2}
{\rm
Consider the system
\begin{equation}\label{lotka-voltera}
\dd x^i(t)=-\alpha_i x^i(t) (1- \langle b, x(t-r)\rangle)\dd t +
\sigma_i x^i(t) (1- \langle b, x(t)\rangle)\, \dd W_i(t),~~ i=1,\ldots, d,
\end{equation}
where $b\in \RR^d_+$, $\alpha_i\ge 0$, $\sigma_i\in\RR$.
The set
$\DD=\left\{ x\in \RR^d_+\, :\,  \langle b, x\rangle\le 1\right\}$
is forward invariant.
Since $\DD$ is a bounded set in $ \RR^d$.
we can modify the nonlinear terms outside some vicinity of $\DD$
in order satisfy the requirement in \eqref{unno}
This allows us to apply Propositions~\ref{pr:wp} and \ref{criteriaGE}(i)
and obtain well-posedness of the problem in \eqref{lotka-voltera}.
The statement on the invariance follows from   Theorem~\ref{th-inv}
via the observation made
in Remark~\ref{re3.3}.
}
\end{Example}

\section{Comparison theorem for sfde's}\label{sect4}
Our next result is  a comparison principle for functional differential equations
perturbed by Kunita type noise of the form \eqref{eq-d}
with the local martingales $M^i$ not only satisfying Hypothesis {\bf (M)} but also
\[
M^i(t,x,\o)= M^i(t,x_i,\o)\quad\mbox{for all}\quad i=1,\ldots,d,\quad x=(x_1,\ldots,x_d),
\]
i.e., $M^i$ depends  on $t$, $\o$  and on the $i$-th
component of spatial variable $x$ only.
As can be seen from \cite{Chu02}, this structural requirement is
needed for a comparison principle even in the non-delay case.
\par
Thus, instead of \eqref{eq-d},
we consider the following Kunita-type retarded stochastic differential
equation
\begin{align}
\label{eq}
\left\{
\begin{array}{rll}
\dd x^i(t)&=&G^i(t,x_t)\,\dd t + M^i(\dd t,x^i(t)),\quad i=1,2,\dots,d, \quad t\ge s,\\
x_s&=&\eta,
\end{array}
\right.
\end{align}
where $\eta$ is a $C$-valued $\F_s$-measurable random variable
and the drift terms $G^i$ satisfies {\bf (G)}.
We fix a decomposition $G=H+b$ as in the previous section and assume that $b^i$ depends on $x^i$ only.
In this case instead of \eqref{psi-d} we have the diagonal system of scalar non-delay equations
\begin{align}
\label{psi}
\left\{
\begin{array}{rll}
\dd \psi^i(t)&=&b^i(t,\psi^i(t))\dd t + M^i(\dd t,\psi^i(t)),\quad i=1,2,\dots,d, \quad t\ge s,\\
\psi^i(s)&=&x\in\RR.
\end{array}
\right.
\end{align}
It follows from Lemma~\ref{le:diff} that
 equation \eqref{psi} generates a stochastic flow of
diffeomorphisms $x \mapsto
\psi^i_{s,t}(x,\o)$  in $\RR$ for each $i=1,2,\dots,d$.
Moreover,
\[
\Psi_{s,t}(x,\o)=(\psi^1_{s,t}(x_1,\o),\ldots,\psi^i_{s,t}(x_d,\o))
\]
satisfies all statements of  Lemma~\ref{le:diff}.
Below we often write $\psi^i(t,x,\o)$ instead of $\psi^i_{0,t}(x,\o)$.
Observe that due to the diffeomorphic property the flow $\psi^i$ is automatically {\em
strongly monotone} in the sense that for  $x,y\in\RR$ we have
$$
x < y \mbox{ implies } \psi^i_{s,t}(x,\o) < \psi^i_{s,t}(y,\o) \mbox{ for each } s,t \ge 0,\,
\o \in \O.
$$
Indeed, if the implication above is not true, then there exist
$x<y$, $s<t$, and $\o$ such that
$ \psi^i_{s,t}(x,\o) = \psi^i_{s,t}(y,\o)$.
Since $\psi^i_{s,t}(\cdot,\o)$ is invertible, this implies $x=y$
and thus provides a contradiction.
\par
Applying Proposition~\ref{pr:equiv-d} we can specify representations
\eqref{int-d1} and \eqref{int-d2} for our  case of diagonal $M^i$.
Namely, if we define
\begin{eqnarray*}
\xi^i(u,x^i,\o)&:=&\psi^i(u,\cdot,\o)^{-1}(x^i)\\
F^i(u,x^i,\eta,\o)&:=&\left\{ \DDD_{x^i}\psi^i(u,x^i,\o) \right\}^{-1} H^i(u,\eta,\o),
\end{eqnarray*}
then \eqref{int-d1} can be written in the form
\begin{equation}
\label{int}
x^i(t,\o)=
\psi^i \left( t,\left[ \xi^i(s,\eta^i(0),\o)+\int_s^t F^i(u,\xi^i (u,x^i(u,\o),\o),x_u(\o),\o)\dd u
\right],\o \right)
\end{equation}
for all $t \ge s$.

Let $C_+$ be the standard cone in $C$. This cone defines a partial
order relation via
\begin{equation}\label{por-C}
\eta\ge\eta_*~~{\rm iff}~~ \eta-\eta_*\in C_+,
\end{equation}
i.e., iff $\eta^i(s)\ge\eta_*^i(s)$ for all $s\in [-r,0]$ and
$i=1,\ldots,d$, where
\[
\eta=(\eta^1,\ldots,\eta^d)
~~~\mbox{and}~~~\eta_*=(\eta_*^1,\ldots,\eta_*^d)
\] are elements
from $C=C([-r,0],\RR^d)$. We write $\eta>\eta_*$ iff
$\eta\ge \eta_*$ and $\eta\ne\eta_*$ and use the notation
$\eta>>\eta_*$ if
\[
\eta^i(s)>\eta_*^i(s)~~{\rm for~all}~~s\in [-r,0]~~
{\rm and}~~i=1,\ldots,d.
\]
We also consider  another sfde
\begin{equation}
\label{eq-cp}
\left\{
\begin{array}{rll}
\dd x^i(t)&=&\bar G^i(t,x_t)\,\dd t + M^i(\dd t,x^i(t)),\quad i=1,2,\dots,d, \quad t\ge s,\\
x_s&=&\eta_*\in C,
\end{array}
\right.
\end{equation}
with the same $M$ and $b$. We assume that the random field
$\bar G=\{\bar G^i\}$ satisfies Hypothesis ({\bf G}) with decomposition $\bar G=\bar H + b$.
Let $\DD (\o)\subseteq\RR^d$  be a closed  set with nonempty interior.
\begin{Definition}
{\rm  Let $\DD (\o)\subseteq\RR^d$  be a closed  set with nonempty interior.
For each $\o$, let $[a(\o),b(\o)]$ a random interval in $\RR_+$.
A random vector field
$G=(G^1,\ldots,G^d)\,:\, [0,\infty) \times C \times \O \to \RR^d$
is said to be {\em quasimonotone} on
$[a(\o),b(\o)]\times \DD \subseteq\RR^{d+1}$
iff for any
$\eta=(\eta^1(s),\ldots,\eta^d(s))$ and
$\eta_*=(\eta_*^1(s),\ldots,\eta_*^d(s))$ from $C_\DD$,
where $C_\DD$ is defined by \eqref{d-0}, we have the following implication
\[
\{\eta\ge \eta_*,\; \eta^i(0)=\eta^i_*(0)\}\;
\Rightarrow\;
G^i(t,\eta,\o)\ge G^i(t,\eta_*,\o)
\]
for every $t\in [a(\o),b(\o)]$, $\o\in\O$ and $i=1,\ldots,d$.
}
\end{Definition}
We note for future use that quasimonotonicity is invariant with respect to a decomposition $G=H+b$ with $b^i$ depending on $x^i$ only
in the sense that $G$ is quasimonotone if and only if  $H$ is quasimonotone.

\begin{Theorem}[Comparison Principle]\label{th-cp1}
Assume that $M$ satisfies the conditions above and that the random vector fields
$G=(G^1,\ldots,G^d)$ and $\bar{G}=(\bar{G}^1,\ldots,\bar{G}^d)$ satisfy Hypothesis {\bf(G)}.
Let $x(t):=x(t,\eta,\o)$ be a solution to \eqref{eq} and
$y(t)=y(t,\eta,\o)$  be a solution to \eqref{eq-cp} which possess the
property
\[
x(t),\; y(t)\in \DD,\quad {for}~~ t\in [s,s+T(\o)]
\]
for some convex  closed set $\DD\subseteq\RR^d$
with nonempty interior, where $T(\o)>0$
for all $\o\in\O$.
Assume that {\bf (M$_\DD$)} holds and the random field $G$ is  quasimonotone on
$[s,s+T(\o)]\times \DD$. Then the following assertions hold:
\begin{enumerate}
\item If $\eta\le\eta_*$ and
\begin{equation}\label{cp-0}
G(t,\xi,\o)\le \bar{G}(t,\xi,\o)~~{for~all}~~
\xi\in C_\DD,\; t\in [s, s+T(\o)],\; \o\in\O,
\end{equation}
then
\begin{equation}\label{cp-1}
x(t;\eta,\o)\le y(t,\eta_*,\o)~{for~all}~~
 t\in [s, s+T(\o)],\; \o\in\O.
\end{equation}
\item If $\eta\ge\eta_*$ and
\begin{equation}\label{cp-2-0}
G(t,\xi,\o)\ge \bar{G}(t,\xi,\o)~~{for~all}~~
\xi\in C_\DD,\; t\in [s, s+T(\o)],\; \o\in\O,
\end{equation}
then
\begin{equation}\label{cp-2}
x(t;\eta,\o)\ge y(t,\eta_*,\o)~{for~all}~~
 t\in [s, s+T(\o)],\; \o\in\O.
\end{equation}
\end{enumerate}
\end{Theorem}

\noindent {\bf Proof.} \hspace{.2cm}
We prove the first part only (the proof of the reversed inequalities
is similar).
\par
We start with the case $\eta<<\eta_*$ and $G<<\bar{G}$, i.e., we assume that
\begin{equation}\label{cp-3}
\eta^i(s)<\eta_*^i(s)~~{\rm for~all}~~s\in [-r,0]~~
{\rm and}~~i=1,\ldots,d,
\end{equation}
and
\begin{equation}\label{cp-4}
G^i(t,\xi,\o)< \bar{G}^i(t,\xi,\o)~~{\rm for~all}~~
\xi\in C_\DD,\; t\in [s, s+T(\o)],\; \o\in\O,
\end{equation}
where $i=1,\ldots,d$. The same is true for $H$ and $\bar{H}$
Let us prove that
\begin{equation}\label{cp-5}
x^i(t;\eta,\o)< y^i(t,\eta_*,\o)~{\rm for~all}~~
 t\in [s, s+T(\o)],\; \o\in\O,\;
i=1,\ldots,d.
\end{equation}
Since $x(t)$ and $y(t)$ are continuous for all $\o\in\O$
relation \eqref{cp-5} is valid for some interval $[s, s+\tau(\o)]$,
where $0<\tau(\o)\le T(\o)$.
If \eqref{cp-5} does not hold for {\em all }
$t$ from  $[s, s+T(\o)]$, then for some $\o$ there exist
$t'\in (0, T(\o))$ and $i\in\{1,\ldots,d\}$ such that
\[
 x^i(t')= y^i(t')~~{\rm and}~~
x^j(t)< y^j(t)~{\rm for~all}~
 t\in [s, s+t'),\;  j=1,\ldots,d.
\]
Using representation \eqref{int} and strict monotonicity of $\psi^j$
we obtain that
\begin{equation}\label{cp-6}
 \zeta_G^i(t')= \zeta_{\bar{G}}^i(t')~~{\rm and}~~
\zeta_G^j(t)< \zeta_{\bar{G}}^j(t)~{\rm for~all}~
 t\in [s, s+t'),\;  j=1,\ldots,d,
\end{equation}
where
\[
\zeta_G^j(t):=\xi^j(s,\eta^j(0),\o)+
\int_s^t F_G^j(u,\xi^j(u,x^j(u),\o),x_u,\o)\dd u,
\]
\[
\zeta_{\bar{G}}^j(t):=\xi^j(s,\eta_*^j(0),\o)+
\int_s^t F_{\bar{G}}^j(u,\xi^j(u,y^j(u),\o),y_u,\o)\dd u
\]
with the following notation:
\begin{eqnarray*}
\xi^j(u,x^j,\o)&:=&\psi^j(u,\cdot,\o)^{-1}(x^j),\\
F_G^j(u,x^j,\eta,\o)&:=&\left\{ \DDD_{x^j}\psi^j(u,x^j,\o) \right\}^{-1}
H^j(u,\eta,\o),
\\
F_{\bar{G}}^j(u,x^j,\eta,\o)&:=&\left\{ \DDD_{x^j}\psi^j(u,x^j,\o) \right\}^{-1}
\bar{H}^j(u,\eta,\o).
\end{eqnarray*}
Since the functions $F_G^j(u,x^j,\eta,\o)$ and $F_{\bar{G}}^j(u,x^j,\eta,\o)$
are continuous for every $\o\in\O$, the  processes
$\zeta_G^j(t)$ and $\zeta_{\bar{G}}^j(t)$ are continuously differentiable
and satisfy the equations
\begin{equation}\label{cp-7}
\frac{\dd}{\dd t}\zeta_G^j(t)= F_G^j(t,\xi^j(t,x^j(t),\o),x_t,\o)
\end{equation}
and
\begin{equation}\label{cp-8}
\frac{\dd}{\dd t}\zeta_{\bar{G}}^j(t)= F_{\bar{G}}^j(t,\xi^j(t,y^j(t),\o),y_t,\o).
\end{equation}
It follows from \eqref{cp-6} that
\[
\zeta_G^i(t')-\zeta_G^i(t)>\zeta_{\bar{G}}^i(t')-\zeta_{\bar{G}}^i(t)
~~{\rm for~all}~~s\le t< s+t'.
\]
This implies that
\begin{equation}\label{cp-9}
\frac{\dd}{\dd t}\zeta_G^i(t')\ge \frac{\dd}{\dd t}\zeta_{\bar{G}}^i(t').
\end{equation}
However, since $x^i(t')=y^i(t')$, from \eqref{cp-8} we have that
\[
\frac{\dd}{\dd t}\zeta_{\bar{G}}^i(t')=
\left\{ \DDD_x\psi^j(u,x^i(t'),\o) \right\}^{-1}H^i(t,y_{t'},\o).
\]
Therefore \eqref{cp-4}  written for $H$ and $\bar{H}$,  quasimonotonicity of
$H$ and \eqref{cp-7}
imply that
\begin{eqnarray*}
\frac{\dd}{\dd t}\zeta_{\bar{G}}^i(t') & >&
\left\{ \DDD_x\psi^j(u,x^i(t'),\o) \right\}^{-1}H^i(t',y_{t'},\o)
\\ &
\ge &
\left\{ \DDD_x\psi^j(u,x^i(t'),\o) \right\}^{-1}H^i(t',x_{t'},\o)=
\frac{\dd}{\dd t}\zeta_G^i(t').
\end{eqnarray*}
This relation  contradicts to \eqref{cp-9}. Thus \eqref{cp-3} and
\eqref{cp-4} imply \eqref{cp-5}.
\par
To prove \eqref{cp-1} for the general case
we first apply the result above to the corresponding equations
with
\[
{\bar{G}}_\eps(t,\eta,\o)={\bar{G}}(t,\eta,\o)+\eps ({\bf e}_1-\eta(0))
\]
and
\[
G_\eps(t,\eta,\o)=G(t,\eta,\o)+\eps ({\bf e}_2-\eta(0)),
\]
where ${\bf e}_1, {\bf e}_2\in \mbox{int}\,\DD$ and ${\bf e}_1<<{\bf e}_2$.
It is clear that \eqref{cp-0} implies that
\[
G^i_\eps(t,\xi,\o)< {\bar{G}}_\eps^i(t,\xi,\o)~~{\rm for~all}~~
\xi\in C_\DD,\; t\in [s, s+T(\o)],\; \o\in\O,\;
i=1,\ldots,d.
\]
Thus, by limit transition we obtain \eqref{cp-1} in the case when
$\eta<<\eta_*$. Using this fact it is easy to prove
 \eqref{cp-1} for every $\eta\le\eta_*$.
\hfill $\Box$

\begin{Remark}\label{re:inv}
{\rm If the  drift term $G$ is quasimonotone on $\RR_+^d$, then
we have that $G(t,\eta,\o)\ge G(t,0,\o)$
for every $\eta\in C_{\RR_+^d}$. Therefore
applying the comparison principle in \eqref{cp-2-0} and \eqref{cp-2}
with $\bar{G}\equiv 0$
we can conclude that
 $\RR_+^d$
is a forward invariant set with respect to
sfde \eqref{eq} when  $G(t,0,\o)\ge 0$ and $M(t,0,\o)\equiv 0$.
}
\end{Remark}

\begin{Example}[Lotka-Volterra type model]\label{ex4.1}
{\rm
Consider the system
\begin{equation}\label{lv}
\left\{
\begin{array}{rcl}
\dd x^i(t) &= &\alpha_i x^i(t) \left(1- \beta_i x^i(t)- \displaystyle\sum_{j=1}^d c_{ij}
\int_{-r}^0 x^j(t+\tau)\dd \mu_{ij}(\tau)\right)\dd t  \\ \\
& & + \sigma_i x^i(t) (R_i- x^i(t)) \dd W_i,~~ t>s,~~ i=1,\ldots, d, \\ [2mm]
x_s &= &\eta\in C.
\end{array}\right.
\end{equation}
Here $\alpha_i$,  $\beta_i$ and $R_i$ are positive numbers,
$c_{ij}\ge 0$, $\sigma_i\in\RR$. We assume that $\mu_{ij}(\tau)$
are left continuous nondecreasing functions on $[-r,0]$ of
bounded variation such that
\[
\mu_{ij}(0)-\mu_{ij}(-r)=1, ~~ i,j=1,\ldots, d.
\]
It is easy to see from Theorem \ref{th-inv} (see also Corollary~\ref{co3.6}
and Remark~\ref{re3.3})  that
$\DD=\prod_{i=1}^d [0,R_i]$ is a
forward invariant set for sfde \eqref{lv} provided $R_i\ge \beta_i^{-1}$
for every $i=1,\ldots, d$.
\par
We note that the global well-posedness of \eqref{lv} follows
from Propositions~\ref{pr:wp} and \ref{criteriaGE}(ii) because we can modify the corresponding drift term outside $\DD$ to satisfy \eqref{unno}.
\par
It is also clear that the functions
\[
\bar{G}^i(\eta):=\alpha_i \eta^i(0) \left(1- \beta_i \eta^i(0)-
\sum_{j=1}^d c_{ij}
\int_{-r}^0 \eta^j(\tau)\dd \mu_{ij}(\tau)\right),
 ~~ i,j=1,\ldots, d,
\]
satisfy the inequality
\[
\alpha_i \eta^i(0) \left(1- \beta_i \eta^i(0)
-\sum_{j=1}^d c_{ij}R_j\right)
\le \bar{G}^i(\eta)\le \alpha_i \eta^i(0) \left(1- \beta_i \eta^i(0)\right)
\]
for every $\eta\in C_\DD$, where $C_\DD$ is given by \eqref{d-0}
with $\DD=\prod_{i=1}^d [0,R_i]$. Since the functions
\[
 G_1^i(\eta):=\alpha_i \eta^i(0) \!\left[1- \beta_i \eta^i(0)
-\!\sum_{j=1}^d c_{ij}R_j\!\right]~\mbox{and}~~
G_2^i(\eta):=\alpha_i \eta^i(0) \left(1- \beta_i \eta^i(0)\right)
\]
are quasimonotone, Theorem~\ref{th-cp1} implies that  for any initial
data $\eta\in C_\DD$ a solution
\[
x(t,\eta,\o)=(x^1(t,\eta,\o),\ldots, x^d(t,\eta,\o))
\]
 to problem \eqref{lv} satisfies the inequality
\begin{equation}\label{compar-lv}
    u^i(t,\eta,\o)\le x^i(t,\eta,\o)\le v^i(t,\eta,\o),
 ~~ i,j=1,\ldots, d,
\end{equation}
where $u(t,\eta,\o)=(u^1(t,\eta,\o),\ldots, u^d(t,\eta,\o))$ solves
the problem
\begin{equation}\label{lv-1}
\left\{
\begin{array}{rcl}
\dd u^i(t) &= &\alpha_i u^i(t) \left(1- \beta_i u^i(t)- \displaystyle \sum_{j=1}^d c_{ij}R_j
\right)\dd t  \\ \\
& & + \sigma_i u^i(t) (R_i- u^i(t)) \dd W_i,~~ t>s,\; i=1,\ldots, d, \\ \\
u^i(0) &= &\min_{\tau\in [-r,0]}\eta^i(s),~~ i=1,\ldots, d,
\end{array}\right.
\end{equation}
and $v(t,\eta,\o)=(v^1(t,\eta,\o),\ldots, v^d(t,\eta,\o))$ solves
the problem
\begin{equation}\label{lv-2}
\left\{
\begin{array}{rcl}
\dd v^i(t) &= &\alpha_i v^i(t) \left(1- \beta_i v^i(t)
\right)\dd t  + \sigma_i v^i(t) (R_i- v^i(t)) \dd W_i,~~ t>s, \\ \\
v^i(0) &= &\max_{\tau\in [-r,0]}\eta^i(s),~~ i=1,\ldots, d.
\end{array}\right.
\end{equation}
We emphasize that problems \eqref{lv-1} and \eqref{lv-2} are direct
sums of one-dimensional ordinary stochastic differential equations.
Long time dynamics of these 1D systems is described with details
(see, e.g., \cite{Chu02} and the references therein).
Thus we can use the relations in \eqref{compar-lv}
to ``localize" dynamics of the original sfde  \eqref{lv}.
}
\end{Example}

\section{Order-preserving RDS generated by sfde's}
\label{sect5}
In this section we consider some other applications of
Theorems~\ref{th-inv} and \ref{th-cp1} from point view of theory of
random dynamical systems (RDS).

\subsection{Generation of RDS in an invariant region}
Following the monograph of Arnold \cite{Arn98}, we introduce the notion of
a random dynamical system.
\begin{Definition}{\rm  Let $X$ be a topological space.
A {\it random dynamical system} (RDS) with time
${\RR}_{+}$ and state space $X$ is a pair $(\vartheta,\phi)$ consisting
of the  following two objects:
\begin{enumerate}
\item A metric dynamical system (MDS)
$\vartheta \equiv (\Omega, {\cal F},{\bf P},\{\vartheta(t), t\in {\RR}\})$,
i.e.,
a probability space $(\Omega, {\cal F},{\bf P})$ with a family of
 measure preserving transformations
$\{\vartheta(t)\, :\, \Omega \mapsto \Omega,  t\in {\RR}\}$ such that
\begin{enumerate}
  \item  $\vartheta(0) =\mathrm{id},\quad \vartheta(t)\circ \vartheta(s) =
\vartheta(t+s)
\quad\mbox{for all}\quad t,s\in {\RR}$;
  \item the map $(t,\omega)\mapsto \vartheta(t)\omega$  is measurable
and $\vartheta(t) {\bf P}={\bf P}$ for all $t \in {\RR}$.
\end{enumerate}
\item A (perfect) cocycle $\phi$ over $\vartheta$ of continuous mappings of $X$
with one-sided time ${\RR}_+$, i.e.
a measurable mapping
$$
\phi\, :\,{\RR}_+\times \Omega\times X\mapsto X,
\quad (t, \omega, x)\mapsto\phi (t, \omega) x
$$
such that
(a)  the mapping  $\phi(\cdot,\omega)\, :\,  x
\mapsto\phi (t, \omega) x$
is continuous for  all $t\ge 0$ and $\omega \in \Omega$;
(b) it satisfies the cocycle property:
$$
\phi(0,\omega) =\mathrm{id},\quad \phi(t+s,\omega)=  \phi(t,
\vartheta(s)\omega)
\circ\phi(s,\omega)
$$
 for all $t,s\ge 0$ and $\omega\in\Omega$.
\end{enumerate}
}
\end{Definition}
\begin{Definition} {\rm
Let $\vth$ be an MDS, $\bar {\cal F}$ the ${\bf P}$-completion
of ${\cal F}$ and ${\bf F}=\{{\cal F}_t ,t \in {\RR}\}$ a family of sub-$\sigma$-algebras of
$\bar {\cal F}$ such that
\begin{enumerate}
\item ${\cal F}_s \subseteq  {\cal F}_t, \quad s < t$;
\item ${\cal F}_s=\bigcap_{h>0}   {\cal F}_{s+h}, \quad s \in \RR\, $,
 i.e. the filtration {\bf F} is right-continuous;
\item ${\cal F}_s$ contains all sets in ${\cal F}$ of ${\bf P}$-measure 0, $s \in {\RR}$;
\item $\vartheta(s)$ is $({\cal F}_{t+s},{\cal F}_t)$-measurable for
all $s,t \in {\RR}$.
\end{enumerate}
Then $(\vth,{\bf F})$ is called a {\em filtered metric dynamical system} (FMDS). If - in
addition - $(\vth,\phi)$ is an RDS such that $\phi(t,\cdot)x$ is
$({\cal F}_t,{\cal B}(X))$-measurable for every $t \ge 0,x \in X$, then
$(\vth,{\bf F},\phi)$ is called a {\em filtered random dynamical
  system} (FRDS).
}
\end{Definition}
We recall that an $X$-valued stochastic process
$Y(t),t\in T \subseteq {\RR}$ is called {\em adapted} or
{\em nonanticipating} with respect to the filtration
${\bf F}$ if $Y(t)$ is $({\cal F}_t,{\cal B}(X))$-
measurable for every $t \in T$.
Therefore $(\vth,{\bf F},\phi)$ is an FRDS iff
$(\vth,\phi)$ is an RDS, $(\vth,{\bf F})$ is an FMDS
and $\phi(\cdot,\cdot)x$ is
adapted to {\bf F} for every $x \in X$.

\begin{Theorem}\label{th-gen}
Assume that Hypotheses {\bf (M$_{\DD}$)} and {\bf (G$_\eps$)} are in force
and   {\bf (GE)} (see Proposition~\ref{criteriaGE})
 holds.
If the drift term $G(t,\eta,\o)\equiv G(\eta)$ and the
the local characteristic $a$ of $M$ are deterministic and autonomous,
then problem \eqref{int-d1},\eqref{int-d2}  (and hence \eqref{eq-d}) generates a FRDS $(\vth, \va)$ in $C_\DD$, where
$C_\DD$ is defined by \eqref{d-0} (the case $\DD\equiv\RR^d$
is not excluded).  The corresponding cocycle $\va$ has the form
\begin{equation*}
[\va(t,\o)\eta](\tau)=
\left\{\begin{array}{cr}
x(t+\tau,\eta,\o), & t+\tau>0, \\
\eta(t+\tau), & t+\tau\le 0,
\end{array}
\right.
\end{equation*}
for every $\tau\in [-r,0]$, where $x(t,\eta,\o)$ is a solution to
problem \eqref{eq-d} with $s=0$. Moreover $\va(t,\o)$ is compact
mapping in $C_\DD$, i.e. for any bounded set $A$ from $C_\DD$,
the set $\va(t,\o)A$ is relatively compact in $C_\DD$ for every $t>0$.
\end{Theorem}
{\bf Proof.}
It follows from Theorem~\ref{th-inv} and from the representation in
 \eqref{int-d1},\eqref{int-d2}
of solutions to \eqref{eq-d}. We also use Propositions~\ref{pr:equiv-d},
\ref{pr:wp} and \ref{criteriaGE}.
\hfill $\Box$
\smallskip\par
To describe long-time dynamics of an RDS we need a notion of a
random set (see, e.g., \cite{Arn98} and the references therein).
\begin{Definition}
{\rm
A mapping $\o\mapsto\ D(\o)$ from $\O$ into the collection of
all subsets of a separable Banach space $V$ is said to be
{\em random closed set}, iff $D(\o)$ is a closed set for any $\o\in\O$
and $\o\mapsto \mbox{dist}_V(x, D(\o))$ is measurable for any
$x\in V$. The random closed set $D(\o)$ is said to be {\em compact}, if
$D(\o)$ is compact for each $\o$.
The random closed set $D(\o)$ is said to be {\em tempered} if
\[
D(\o)\subset\{ x\in V\, :\, \| x\|_V\le r(\o)\},~~ \o\in\O,
\]
where the random variable $r(\o)$ possesses the property
$\sup_{t\in \RR}\{ r(\vth(t)\o)e^{-\gamma |t|}\}<\infty$
for any $\gamma>0$.
}
\end{Definition}

We also need the following concept of a random attractor of an RDS
(see \cite{CF,Sch1} and also \cite{Arn98,Chu02} and the
references therein).
Below we denote by $X$ any subset af separable Banach space $V$
equipped with the induced topology.
\par
Let $\sD$ be a family of
random closed sets in $X$ which is closed with respect to inclusions (i.e. if
$D_1\in \sD$ and a random closed set $\{D_2(\o)\}$ possesses the property
$D_2(\omega)\subset D_1(\omega)$ for all $\o\in\O$, then
$D_2\in \sD$). Sometimes the collection $\sD$ is called a {\em universe}
of sets (see \cite{Arn98}).

\begin{Definition}\label{d1.6.1}
{\rm
Suppose that $(\vth,\va)$ is an RDS in $X$.
 Let $\sD$ be a universe.
A random closed set $\{ A(\o)\}$ from $\sD$
is said to be a {\em random pull-back attractor} of the  RDS
$(\vth,\va)$
in $\sD$ if $A(\o)\neq X$ for every $\o\in\O$ and
the following properties hold:
\begin{enumerate}
\item[(i)] $A$ is an invariant set, i.e.
 $\varphi(t,\omega)A(\omega) =A(\vartheta(t)\omega)$ for $t\ge 0$ and
$\omega\in\Omega$;
\item[(ii)] $A$ is attracting in $\sD$, i.e. for all $D\in \sD$
\begin{equation}\label{1.6.1}
\lim_{t\to +\infty}
d_X\{\varphi(t,\vartheta(-t)\omega)D(\vartheta(-t)\omega)\, | \,
A(\omega)\}=0, \quad \o\in\O\;,
\end{equation}
where $d_X\{A | B\}= \sup_{x\in A} {\rm dist}_X (x, B)$.
\end{enumerate}
If instead of \eqref{1.6.1} we have that
\begin{equation*}
\lim_{t\to +\infty}\PP\left\{ \o\,:\;
d_X\{\varphi(t,\omega)D(\o)\, |\,
A(\vth(t)\omega)\}\ge \delta\right\}=0
\end{equation*}
for any $\delta>0$,   then
is said to be a {\em random weak attractor} of the  RDS
$(\vth,\va)$.
}
\end{Definition}
Some authors (e.g.~\cite{CF,Arn98}) require a random attractor to be compact (and do not insist
that it is different from the whole space). This distinction will not be important in what follows.
The notion of a  {\em weak random attractor} was introduced in \cite{Ochs}.
For the relation between weak, pull-back and forward attractors
we refer to \cite{Scheu2}.

\begin{Remark}\label{re:attractor}
{\rm
If $\DD$ is bounded in $\RR^d$, then it is easy to see that
the  RDS $(\vth,\va)$ generated by \eqref{eq-d}
in $C_\DD$
has  a  random compact pull-back attractor
in the universe $\sD$ all  bounded sets.
Since by  Theorem~\ref{th-gen} the   RDS $(\vth,\va)$
is compact,
this follows from Theorem~1.8.1 \cite{Chu02}, for instance.
In the case of unbounded sets $\DD$ (e.g., $\DD=\RR^d_+$)
we need some conditions which guarantee dissipativity
of the corresponding RDS. These conditions can be obtained in the same
way as for the non-delay case (see, for instance,
 \cite[Theorem 6.5.1]{Chu02}).
}
\end{Remark}

\subsection{Monotone RDS}
Let as above $C=C([-r,0],\RR^d)$ and
 $C_+$ be the standard cone in $C$
of nonnegative elements:
\[
C_+=\left\{ \eta=(\eta^1(s),\ldots,\eta^d(s))\in C~~
\eta^i(s)\ge 0~~ \forall~s\in [-r,0],~~i=1,\ldots,d
\right\}.
\]
This cone  is a normal solid minihedral cone.
This fact is important for further application of the theory
of monotone RDS.
We refer to \cite{Kr1} and \cite{Kr4} for more details
concerning cones and partially ordered  spaces.
\par
Let $\DD$ be a convex closed set in $\RR^d$
with nonempty interior. In
the space $C_\DD$  given by \eqref{d-0}
we we define a partial
order relation via \eqref{por-C}, i.e.,
$\eta\ge\eta_*$  iff $\eta-\eta_*\in C_+$.

\begin{Theorem}\label{th-mrs}
Assume that the hypotheses of Theorem~\ref{th-gen} hold.
Let $\DD$ be a convex closed set in $\RR^d$
with nonempty interior and $(\vth,\va)$ be  the FRDS generated
by problem \eqref{eq} in $C_\DD$ defined by \eqref{d-0}.
If the random field $G$ is  quasimonotone in $\DD$
then $(\vth,\va)$  is an order-preserving FRDS in $C_\DD$  which
means that
$$
\eta\le \xi~\mbox{in}~C_\DD  \quad\mbox{implies}\quad \va (t,\omega)\eta\le\va (t,\omega)\xi
\quad\mbox{for all}~ t\ge 0~\mbox{and}~ \omega\in\Omega.
$$
\end{Theorem}
{\bf Proof.} This follows from Theorem~\ref{th-cp1} with $\bar G\equiv G$.
\hfill $\Box$

Theorem~\ref{th-mrs} makes it possible  to apply the general theory
of monotone RDS (see \cite{Chu02} and also \cite{ArChu98a,ChuSch03a})
to the class of sfde's considered.
In particular it is possible to obtain  the following results:
\begin{itemize}
  \item To provide transparent conditions which guarantee
  the existence of sto\-chastic equilibria and a compact pull-back attractor
  (see, e.g., the general  Theorem~3.5.1 in \cite{Chu02}).
  We recall (see \cite{Arn98}) that   a random variable $u:\Omega\mapsto C_\DD$ is said
to be  an {\em equilibrium} (or fixed point, or stationary solution)
of the RDS $(\vartheta,\va)$ if it is invariant under $\va$, i.e.\ if
\[
\va (t, \omega)u(\omega)=u(\vartheta(t)\omega)
\quad\mbox{a.s. for all}\quad t\ge0.
\]
One can see that any equilibrium $u(\o)\in C_\DD$ for  $(\vth,\va)$
has the form
$u(\tau,\o)= v(\vth(\tau)\o)$, $\tau\in [-r,0]$,
where $v(\o)$ is a random variable in $\DD$.
  \item To describe the pull-back attractor for $(\vartheta,\va)$ as a compact set
  lying between  two of  its equilibria (see general Theorem~3.6.2 in \cite{Chu02}).
  \item  In the case when
there exists a probability measure $\pi$ on
 the Borel $\sigma$-algebra of subsets $C_\DD$,
such that the law ${\cal L}(\va(t,\o)x)$ weakly converges to
$\pi$  in $C_\DD$, the system
   $(\vth,\va)$ has a random weak attractor $A(\o)$ which  is singleton, i.e. $A(\o)=\{v(\o)\}$,
where the random variable $v(\o)\in C_\DD$ is an equilibrium
(see the general Theorem 1 proved in \cite{ChuSch03a}).
We  note that
sufficient conditions for the existence of an invariant probability measure $\pi$
and weak convergence of transition probabilities to $\pi$ for an sfde (monotone or not) have been established
for example in \cite {IN64}, \cite{Scheu0}, \cite{ESvG10}, and \cite{HMS11}.
\end{itemize}

\begin{Example}\label{ex5.1}
{\rm
Consider the stochastic equations
\begin{equation}\label{5.3}
\dd x^i(t)=
(g^i_0(x^1(t),\ldots,x^d(t))+g^i_1(x^1(t-r_1),\ldots,x^d(t-r_d))\dd t
+m^i(x^i)\dd W_i
\end{equation}
for $ i=1,\ldots, d$.
We assume that $g_0^i$, $g_1^i$ and  $m^i$
are smooth functions which are globally Lipschitz. Under these conditions
we can apply Proposition~\ref{criteriaGE}(iii) to guarantee
global well-posedness for \eqref{5.3}.
Moreover, one can see from Theorem~\ref{th-mrs} that
 equations \eqref{5.3} generate an order-preserving RDS in the space $C=C([-r,0];\RR^d)$
with $r=\max_i r_i$ provided that
\begin{equation*}
\frac{\partial g_0^i(x)}{\partial x_i}\ge 0,~~ x\in\RR^d,~
i\neq j,
\end{equation*}
and $g_1^i(x)$ is monotone, i.e. for every $i=1,\ldots, d$ we have
that
\begin{equation*}
g^i_1(x^1,\ldots,x^d)\le g^i_1(y^1,\ldots, y^d)~~\mbox{when}
~~x^j\le y^j,\; j=1,\ldots, d.
\end{equation*}
}
\end{Example}

The following example is a special case of Example~\ref{ex5.1}.

\begin{Example}\label{ex5.1-1D}
{\rm
Let $W$ be standard Brownian motion. Consider 1D the retarded stochastic differential equation
\[
\dd x(t)=(f(x(t))+g(x(t-1)))\dd t+ \sigma(x(t))\dd W(t),
\]
where $f,g,\sigma:\RR \to \RR$ are Lipschitz, $g$ is monotone and
$\sigma$ is strictly positive.
By Theorem~\ref{th-mrs}
this equation generates an order-preserving RDS in $C([-1,0],\RR)$
Assume that the associated Markov semigroup on $C([-1,0],\RR)$ admits an invariant (or {\em stationary})
measure (sufficient conditions are provided in \cite{Scheu0}).
In this case we can apply Theorem 1 \cite{ChuSch03a}
and conclude that
the corresponding RDS  has a unique equilibrium which is a weak random
attractor.

It is known (cf. \cite{KWW99}) that
 in the case $\sigma\equiv 0$ the attractor of this system
can contain multiple equilibria and also a periodic orbit.
Thus
we observe here that adding the noise term simplifies essentially
long-time behavour  of the system (see also \cite{ChuSch03a} for some details).
}
\end{Example}

\begin{Example}[Stochastic biochemical control circuit]\label{ex6.1}
{\rm We
consider the following system of   Stratonovich  stochastic
equations
\begin{equation}\label{6.10}
\dd x^1(t) =(g(L_d x^d_t)-\alpha_1 x^1(t))\dd t+ \sigma_1\cdot
  x^1(t)\circ \dd W_t^1\;,
\end{equation}
\begin{equation}\label{6.11}
\dd x^j(t) =(L_{j-1}x^{j-1}_t-\alpha_j x^j(t))\dd t+ \sigma_j\cdot
x^j(t)\circ \dd W_t^j,\quad j=2,\ldots, d\;.
\end{equation}
Here as above ``$\circ$'' denotes Stratonovich integration,
 $\sigma_j$ are nonnegative and $\alpha_j$ are
positive constants, $j=1,\ldots, d$, and $g : {\RR_+} \mapsto {\RR_+}$ is a
$C^1$ function such that
\begin{equation*}
0<  g(u)\le au+b,\quad\mbox{and}\quad  g'(u)\ge 0 \quad\mbox{for
every}\quad u>0
\end{equation*}
for some constants $a$ and $b$. We also use
the notation $x_{j,t}(s)=x_j(t+s)$ for $s\in [-r_j,0]$
and
\[
L_j\eta=\int_{-r_j}^0\eta(s) \dd \mu_j(s),
\]
where $\mu_j : [-r_j,0]\mapsto\RR$ is nondecreasing,
$\mu_j(-r_j)=0$, $\mu_j(0)=1$, $\mu_j(s)>0$ for
$s>-r_j$.
We   denote $r=\max_jr_j$ and  equip \eqref{6.10} and \eqref{6.11} with initial data
\begin{equation}\label{6-id}
    x^i(t)=\xi^i(t)\ge 0,~~~t\in [-r,0],~~~ i=1,\ldots,d.
\end{equation}
A deterministic version of this system was considered in \cite{Sm},
the stochastic non-retarded case was studied in \cite{Chu02},
see also \cite{ChuSch03a}.
\par
Let $\xi\in C([-r,0]:\RR^d_+)$ by Proposition ~\ref{pr:wp} a local solution
$$x(t)=(x^1(t),\ldots,x^d(t))$$ exists on some interval $[0,T(\o))$.
By Theorem~\ref{th-inv} (see also Remark~\ref{re3.2} and Corollary~\ref{co3.6})
we have that $\RR_+^d$ is a forward invariant set, i.e.
$x(t)\in \RR_+^d$ for all $t\in [0,T(\o))$.
Applying Comparison Principle (see Theorem~\ref{th-cp1})
we conclude that
\begin{equation}\label{cp-bio}
    0\le x(t)\le \bar{x}(t)~~~\mbox{for all}~~t\in [0,T(\o)),
\end{equation}
where $\bar{x}^i(t)=(\bar{x}^1(t),\ldots,\bar{x}^d(t))$
solves
the following system of linear equations
\begin{equation}\label{6.10a}
\dd x^1(t) =(a L_d x^d_t-\alpha_1 x^1(t)+ b)\dd t+ \sigma_1\cdot
  x^1(t)\circ \dd W_t^1\;,
\end{equation}
\begin{equation}\label{6.11a}
\dd x^j(t) =(L_{j-1}x^{j-1}_t-\alpha_j x^j(t))\dd t+ \sigma_j\cdot
x^j(t)\circ \dd W_t^j,\quad j=2,\ldots, d\;.
\end{equation}
with initial data \eqref{6-id}. The structure of \eqref{6.10a}
and \eqref{6.11a} allows us to solve these equations.
Indeed, if we  consider the drift part of the problem:
\[
\dd x^j(t) = \sigma_j\cdot
x^j(t)\circ \dd W_t^j,~~~ x^j(0)=x,\quad j=1,\ldots, d\;.,
\]
then $\psi^j(t,x)=x \exp\{\sigma_j W^j(t)\}$
 solves it. The spatial derivative of $\psi^j(t,x)$ and its inverse are
 are independent of $x$ and thus we can apply
 Proposition~\ref{criteriaGE}(iii) to prove
 global existence of the solution  $\bar{x}$.
Due to \eqref{cp-bio} this implies that the solution $x(t)$
of  \eqref{6.10}, \eqref{6.11} and \eqref{6-id} does not explode
at finite time and thus
equations \eqref{6.10} and \eqref{6.11} generate an RDS $(\vth,\va)$ in $C_+=C([-r,0]:\RR^d_+)$, where $r=\max_i r_i$.
By Theorem~\ref{th-mrs} this RDS is order-preserving.
By  Comparison Theorem~\ref{th-cp1}, this system is dominated from
above by the affine RDS
$(\vth,\va_{af})$ generated by \eqref{6.10a}
and \eqref{6.11a}.
\par
Now we concentrate  on
 the case $a=0$ (this means that $g(u)$ is bounded). In this case  we can construct
an equilibrium $v(\o)=(v^1(\o),\ldots, v^d(\o))$ for $(\vth,\va_{af})$
by the formulas
\[
v_1(\o)=b\int^0_{-\infty}e^{\alpha_1 t-\sigma_1W^1_t} \dd t,
\]
and
\[
v_j(\o)=\int^0_{-\infty}L_{j-1}v^{j-1}_t\cdot e^{\alpha_j t-\sigma_jW^j_t} \dd t,
~~j=2,\ldots, d,
\]
where $v_t^j(\o):= v^j(\vth(t+\tau))$, $\tau\in [-r,0]$.
Since $\va(t,\o)x\le\va_{af}(t,\o)x$ for every $x\in C_+$,
it is easy to see that $v(\o)$ is a super-equilibrium for
$(\vth,\va)$, i.e.,
\[
\va (t, \omega)u(\omega)\le v(\vartheta_t\omega)
\quad\mbox{a.s. for all}\quad t\ge0.
\]
Thus by Theorem 3.5.1~\cite{Chu02} the RDS $(\vth,\va)$
has an equilibrium $u(\o)\in \RR^d_+$. If $g(0)>0$, this equilibrium
is strongly positive.
\par
We can also show that in the case $a=0$ the RDS $(\vth,\va)$
possesses a random pull-back attractor in the universe of
all tempered subsets of $C([-r,0]:\RR^d_+)$.
Indeed, due to the compactness
property of the cocycle $\vth$ (see Theorem~\ref{th-gen}) it is sufficient
to prove that $(\vth,\va)$ possesses a bounded a
absorbing set. This set can be constructed
in the following way.
\par
Let $v(\o)$ be the equilibrium for $(\vth,\va_{af})$ constructed above.
One can see that in this case $v_\la(\o)=\la v(\o)$ is
a super-equilibrium for RDS $(\vth,\va_{af})$ for every $\la>1$.
One can also see that the top Lyapunov exponent
for $(\vth,\va_{af})$ with $a=b=0$ is negative. This implies
$v_\la(\o)$ is an {\em  absorbing}   super-equilibrium  for $(\vth,\va_{af})$, i.e.
for every tempered set $D(\o)$ in $C([-r,0]:\RR^d_+)$ there is
$t_D(\o)$ that
\[
\va_{af}(t,\vth_{-t}\o)y(\vth_{-t}\o)\le v_\la(\o),~~~ t\ge t_D(\o),~~ y\in D.
\]
Since   $(\vth,\va_{af})$ dominates $(\vth,\va)$, this implies that
the interval
\[
[0, v_\la(\o)]=\left\{ u\in C\, :\; 0\le u\le v_\la(\o)\right\}
\]
 is absorbing for $(\vth,\va)$.
Therefore Theorem~1.8.1\cite{Chu02} implies the existence
of a pullback attractor which belongs to
some interval of the form $[u_1(\o), u_2(\o)]$,
where $[u_1(\o)$ and $u_2(\o)]$ are two equilibria such that
$0\le u_1(\o)\le u_2(\o)\le v(\o)$.
}
\end{Example}

\end{document}